%%%%%%%%%%%%%%%%%%  tex macros for preprints, cm version %%%%%%%%%%%%%%
%                     (P. Ginsparg, last updated 9/91)
%                if confused, type `b' in response to query
%
%---------------------------------------------------------------------%
%% site dependent options:
%% \unredoffs and \redoffs define horizontal and vertical offsets
%% respectively for unreduced and reduced modes. \speclscape defines
%% the \special{} call that sets printer to landscape (sideways) mode.
%% from standard set below, leave uncommented as appropriate or redefine
%
%%% next 400dpi
%\def\unredoffs{} \def\redoffs{\voffset=-.31truein\hoffset=-.48truein}
%\def\speclscape{\special{landscape}}
%
%%% apple lw
\def\unredoffs{} \def\redoffs{\voffset=-.31truein\hoffset=-.59truein}
\def\speclscape{\special{ps: landscape}}
%
%%% qms lasergrafix:
%\def\unredoffs{} \def\redoffs{\voffset=-.4truein\hoffset=.125truein}
%\def\speclscape{\special{qms: landscape}}
%
%%% saclay A4 paper:
%\def\unredoffs{\hoffset-.14truein\voffset-.2truein}
%\def\redoffs{\voffset=-.45truein\hoffset=-.21truein}
%\def\speclscape{\special{landscape}}
%
%---------------------------------------------------------------------%
%
\newbox\leftpage \newdimen\fullhsize \newdimen\hstitle \newdimen\hsbody
\tolerance=1000\hfuzz=2pt
\catcode`\@=11 % This allows us to modify PLAIN macros.
%\def\bigans{b }    %comment two lines to remove option
%\message{ big or little (b/l)? }\read-1 to\answ
%
\ifx\answ\bigans\message{(This will come out unreduced.}
\magnification=1200\unredoffs\baselineskip=16pt plus 2pt minus 1pt
\hsbody=\hsize \hstitle=\hsize %take default values for unreduced format
\else\message{(This will be reduced.} \let\l@r=L
\magnification=1000\baselineskip=16pt plus 2pt minus 1pt
\vsize=7truein \redoffs
\hstitle=8truein\hsbody=4.75truein\fullhsize=10truein\hsize=\hsbody
\output={\ifnum\pageno=0 %%% This is the HUTP version
  \shipout\vbox{\speclscape{\hsize\fullhsize\makeheadline}
    \hbox to \fullhsize{\hfill\pagebody\hfill}}\advancepageno
  \else
  \almostshipout{\leftline{\vbox{\pagebody\makefootline}}}\advancepageno
  \fi}
\def\almostshipout#1{\if L\l@r \count1=1 \message{[\the\count0.\the\count1]}
      \global\setbox\leftpage=#1 \global\let\l@r=R
 \else \count1=2
  \shipout\vbox{\speclscape{\hsize\fullhsize\makeheadline}
      \hbox to\fullhsize{\box\leftpage\hfil#1}}  \global\let\l@r=L\fi}
\fi
%---------------------------------------------------------------------
%
\newcount\yearltd\yearltd=\year\advance\yearltd by -1900

\def\Title#1#2{\nopagenumbers\abstractfont\hsize=\hstitle\rightline{#1}%
\vskip 1in\centerline{\titlefont #2}\abstractfont\vskip
.5in\pageno=0}
\def\Date#1{\vfill\leftline{#1}\tenpoint\supereject\global\hsize=\hsbody%
\footline={\hss\tenrm\folio\hss}}%  restores pagenumbers
%
%       use following instead of \Date on the preliminary draft,
%       puts date/time on each page in big mode, writes labels in margins

\def\draftmode{\message{ DRAFTMODE }\def\draftdate{{\rm preliminary draft:
\number\month/\number\day/\number\yearltd\ \ \hourmin}}%
\headline={\hfil\draftdate}\writelabels\baselineskip=20pt plus 2pt
minus 2pt
 {\count255=\time\divide\count255 by 60 \xdef\hourmin{\number\count255}
  \multiply\count255 by-60\advance\count255 by\time
  \xdef\hourmin{\hourmin:\ifnum\count255<10 0\fi\the\count255}}}
%       use \nolabels to get rid of eqn, ref, and fig labels in draft mode
\def\nolabels{\def\wrlabeL##1{}\def\eqlabeL##1{}\def\reflabeL##1{}}
\def\writelabels{\def\wrlabeL##1{\leavevmode\vadjust{\rlap{\smash%
{\line{{\escapechar=` \hfill\rlap{\sevenrm\hskip.03in\string##1}}}}}}}%
\def\eqlabeL##1{{\escapechar-1\rlap{\sevenrm\hskip.05in\string##1}}}%
\def\reflabeL##1{\noexpand\llap{\noexpand\sevenrm\string\string\string##1}}}
\nolabels
%
% tagged sec numbers
\global\newcount\secno \global\secno=0 \global\newcount\meqno
\global\meqno=1
\def\newsec#1{\global\advance\secno by1\message{(\the\secno. #1)}
%\ifx\answ\bigans \vfill\eject \else \bigbreak\bigskip \fi  %if desired
\global\subsecno=0\eqnres@t\noindent{\bf\the\secno. #1}
\writetoca{{\secsym} {#1}}\par\nobreak\medskip\nobreak}
\def\eqnres@t{\xdef\secsym{\the\secno.}\global\meqno=1\bigbreak\bigskip}
\def\sequentialequations{\def\eqnres@t{\bigbreak}}\xdef\secsym{}
\global\newcount\subsecno \global\subsecno=0
\def\subsec#1{\global\advance\subsecno by1\message{(\secsym\the\subsecno. #1)}
\ifnum\lastpenalty>9000\else\bigbreak\fi
\noindent{\it\secsym\the\subsecno. #1}\writetoca{\string\quad
{\secsym\the\subsecno.} {#1}}\par\nobreak\medskip\nobreak}
\def\appendix#1#2{\global\meqno=1\global\subsecno=0\xdef\secsym{\hbox{#1.}}
\bigbreak\bigskip\noindent{\bf Appendix #1. #2}\message{(#1. #2)}
\writetoca{Appendix {#1.} {#2}}\par\nobreak\medskip\nobreak}
%
%       \eqn\label{a+b=c}   gives displayed equation, numbered
%               consecutively within sections.
%     \eqnn and \eqna define labels in advance (of eqalign?)
%
\def\eqnn#1{\xdef #1{(\secsym\the\meqno)}\writedef{#1\leftbracket#1}%
\global\advance\meqno by1\wrlabeL#1}
\def\eqna#1{\xdef #1##1{\hbox{$(\secsym\the\meqno##1)$}}
\writedef{#1\numbersign1\leftbracket#1{\numbersign1}}%
\global\advance\meqno by1\wrlabeL{#1$\{\}$}}
\def\eqn#1#2{\xdef #1{(\secsym\the\meqno)}\writedef{#1\leftbracket#1}%
\global\advance\meqno by1$$#2\eqno#1\eqlabeL#1$$}
%
%            footnotes
\newskip\footskip\footskip14pt plus 1pt minus 1pt %sets footnote baselineskip
\def\footnotefont{\ninepoint}\def\f@t#1{\footnotefont #1\@foot}
\def\f@@t{\baselineskip\footskip\bgroup\footnotefont\aftergroup\@foot\let\next}
\setbox\strutbox=\hbox{\vrule height9.5pt depth4.5pt width0pt}
\global\newcount\ftno \global\ftno=0
\def\foot{\global\advance\ftno by1\footnote{$^{\the\ftno}$}}
%
%say \footend to put footnotes at end
%will cause problems if \ref used inside \foot, instead use \nref before
\newwrite\ftfile
\def\footend{\def\foot{\global\advance\ftno by1\chardef\wfile=\ftfile
$^{\the\ftno}$\ifnum\ftno=1\immediate\openout\ftfile=foots.tmp\fi%
\immediate\write\ftfile{\noexpand\smallskip%
\noexpand\item{f\the\ftno:\ }\pctsign}\findarg}%
\def\footatend{\vfill\eject\immediate\closeout\ftfile{\parindent=20pt
\centerline{\bf Footnotes}\nobreak\bigskip\input foots.tmp }}}
\def\footatend{}
%
%     \ref\label{text}
% generates a number, assigns it to \label, generates an entry.
% To list the refs on a separate page,  \listrefs
%
\global\newcount\refno \global\refno=1
\newwrite\rfile
%
%\def\ref{[\the\refno]\nref}  %This line is replaced by
                              %the next line when references
                              %are put alphabetically at the
                              %beginning of the paper
\def\ref{\nref}
\def\nref#1{\xdef#1{[\the\refno]}\writedef{#1\leftbracket#1}%
\ifnum\refno=1\immediate\openout\rfile=refs.tmp\fi
\global\advance\refno by1\chardef\wfile=\rfile\immediate
\write\rfile{\noexpand\item{#1\
}\reflabeL{#1\hskip.31in}\pctsign}\findarg}
%   horrible hack to sidestep tex \write limitation
\def\findarg#1#{\begingroup\obeylines\newlinechar=`\^^M\pass@rg}
{\obeylines\gdef\pass@rg#1{\writ@line\relax #1^^M\hbox{}^^M}%
\gdef\writ@line#1^^M{\expandafter\toks0\expandafter{\striprel@x #1}%
\edef\next{\the\toks0}\ifx\next\em@rk\let\next=\endgroup\else\ifx\next\empty%
\else\immediate\write\wfile{\the\toks0}\fi\let\next=\writ@line\fi\next\relax}}
\def\striprel@x#1{} \def\em@rk{\hbox{}}
\def\lref{\begingroup\obeylines\lr@f}
\def\lr@f#1#2{\gdef#1{\ref#1{#2}}\endgroup\unskip}

\def\addref#1{\immediate\write\rfile{\noexpand\item{}#1}} %now unnecessary
\def\startrefs#1{\immediate\openout\rfile=refs.tmp\refno=#1}
\def\refs#1{\count255=1[\r@fs #1{\hbox{}}]}
\def\r@fs#1{\ifx\und@fined#1\message{reflabel \string#1 is undefined.}%
\nref#1{need to supply reference \string#1.}\fi%
\vphantom{\hphantom{#1}}\edef\next{#1}\ifx\next\em@rk\def\next{}%
\else\ifx\next#1\ifodd\count255\relax\xref#1\count255=0\fi%
\else#1\count255=1\fi\let\next=\r@fs\fi\next}
%

%
% this is ugly, but moore insists
\newwrite\ffile\global\newcount\figno \global\figno=1
\def\fig{fig.~\the\figno\nfig}
\def\nfig#1{\xdef#1{fig.~\the\figno}%
\writedef{#1\leftbracket fig.\noexpand~\the\figno}%
\ifnum\figno=1\immediate\openout\ffile=figs.tmp\fi\chardef\wfile=\ffile%
\immediate\write\ffile{\noexpand\medskip\noexpand\item{Fig.\
\the\figno. }
\reflabeL{#1\hskip.55in}\pctsign}\global\advance\figno
by1\findarg}
\def\xfig{\expandafter\xf@g}\def\xf@g fig.\penalty\@M\ {}
\def\figs#1{figs.~\f@gs #1{\hbox{}}}
\def\f@gs#1{\edef\next{#1}\ifx\next\em@rk\def\next{}\else
\ifx\next#1\xfig #1\else#1\fi\let\next=\f@gs\fi\next}
\newwrite\lfile
{\escapechar-1\xdef\pctsign{\string\%}\xdef\leftbracket{\string\{}
\xdef\rightbracket{\string\}}\xdef\numbersign{\string\#}}

\def\writestop{\def\writestoppt{\immediate\write\lfile{\string\pageno%
\the\pageno\string\startrefs\leftbracket\the\refno\rightbracket%
\string\def\string\secsym\leftbracket\secsym\rightbracket%
\string\secno\the\secno\string\meqno\the\meqno}\immediate\closeout\lfile}}
\def\writestoppt{}\def\writedef#1{}
\def\seclab#1{\xdef #1{\the\secno}\writedef{#1\leftbracket#1}\wrlabeL{#1=#1}}
\def\subseclab#1{\xdef #1{\secsym\the\subsecno}%
\writedef{#1\leftbracket#1}\wrlabeL{#1=#1}}
\newwrite\tfile \def\writetoca#1{}
\def\leaderfill{\leaders\hbox to 1em{\hss.\hss}\hfill}
%   use this to write file with table of contents
\def\writetoc{\immediate\openout\tfile=toc.tmp
   \def\writetoca##1{{\edef\next{\write\tfile{\noindent ##1
   \string\leaderfill {\noexpand\number\pageno} \par}}\next}}}
%       and this lists table of contents on second pass
%\def\listtoc{\centerline{\bf Contents}\nobreak\medskip{\baselineskip=12pt
%
%  1/98 "toc.tex" has been replaced by "toc.tmp". To get
%  table of content, insert the lines \listtoc and \writetoc, tex twice.
%
% \parskip=0pt\catcode`\@=11 \input toc.tex \catcode`\@=12 \bigbreak\bigskip}}
%

\catcode`\@=12 % at signs are no longer letters
%
%   Unpleasantness in calling in abstract and title fonts
\edef\tfontsize{\ifx\answ\bigans scaled\magstep3\else
scaled\magstep4\fi} \font\titlerm=cmr10 \tfontsize
\font\titlerms=cmr7 \tfontsize \font\titlermss=cmr5 \tfontsize
\font\titlei=cmmi10 \tfontsize \font\titleis=cmmi7 \tfontsize
\font\titleiss=cmmi5 \tfontsize \font\titlesy=cmsy10 \tfontsize
\font\titlesys=cmsy7 \tfontsize \font\titlesyss=cmsy5 \tfontsize
\font\titleit=cmti10 \tfontsize \skewchar\titlei='177
\skewchar\titleis='177 \skewchar\titleiss='177
\skewchar\titlesy='60 \skewchar\titlesys='60
\skewchar\titlesyss='60
\def\titlefont{\def\rm{\fam0\titlerm}% switch to title font
\textfont0=\titlerm \scriptfont0=\titlerms
\scriptscriptfont0=\titlermss \textfont1=\titlei
\scriptfont1=\titleis \scriptscriptfont1=\titleiss
\textfont2=\titlesy \scriptfont2=\titlesys
\scriptscriptfont2=\titlesyss \textfont\itfam=\titleit
\def\it{\fam\itfam\titleit}\rm}
 \ifx\answ\bigans\else scaled\magstep1\fi
\ifx\answ\bigans\def\abstractfont{\tenpoint}\else
\font\abssl=cmsl10 scaled \magstep1 \font\absrm=cmr10
scaled\magstep1 \font\absrms=cmr7 scaled\magstep1
\font\absrmss=cmr5 scaled\magstep1 \font\absi=cmmi10
scaled\magstep1 \font\absis=cmmi7 scaled\magstep1
\font\absiss=cmmi5 scaled\magstep1 \font\abssy=cmsy10
scaled\magstep1 \font\abssys=cmsy7 scaled\magstep1
\font\abssyss=cmsy5 scaled\magstep1 \font\absbf=cmbx10
scaled\magstep1 \skewchar\absi='177 \skewchar\absis='177
\skewchar\absiss='177 \skewchar\abssy='60 \skewchar\abssys='60
\skewchar\abssyss='60
\def\abstractfont{\def\rm{\fam0\absrm}% switch to abstract font
\textfont0=\absrm \scriptfont0=\absrms \scriptscriptfont0=\absrmss
\textfont1=\absi \scriptfont1=\absis \scriptscriptfont1=\absiss
\textfont2=\abssy \scriptfont2=\abssys \scriptscriptfont2=\abssyss
\textfont\itfam=\bigit \def\it{\fam\itfam\bigit}\def\footnotefont{\tenpoint}%
\textfont\slfam=\abssl \def\sl{\fam\slfam\abssl}%
\textfont\bffam=\absbf \def\bf{\fam\bffam\absbf}\rm}\fi
\def\tenpoint{\def\rm{\fam0\tenrm}% switch back to 10-point type
\textfont0=\tenrm \scriptfont0=\sevenrm \scriptscriptfont0=\fiverm
\textfont1=\teni  \scriptfont1=\seveni  \scriptscriptfont1=\fivei
\textfont2=\tensy \scriptfont2=\sevensy \scriptscriptfont2=\fivesy
\textfont\itfam=\tenit \def\it{\fam\itfam\tenit}\def\footnotefont{\ninepoint}%
\textfont\bffam=\tenbf
\def\bf{\fam\bffam\tenbf}\def\sl{\fam\slfam\tensl}\rm}
\font\ninerm=cmr9 \font\sixrm=cmr6 \font\ninei=cmmi9
\font\sixi=cmmi6 \font\ninesy=cmsy9 \font\sixsy=cmsy6
\font\ninebf=cmbx9 \font\nineit=cmti9 \font\ninesl=cmsl9
\skewchar\ninei='177 \skewchar\sixi='177 \skewchar\ninesy='60
\skewchar\sixsy='60
\def\ninepoint{\def\rm{\fam0\ninerm}% switch to footnote font
\textfont0=\ninerm \scriptfont0=\sixrm \scriptscriptfont0=\fiverm
\textfont1=\ninei \scriptfont1=\sixi \scriptscriptfont1=\fivei
\textfont2=\ninesy \scriptfont2=\sixsy \scriptscriptfont2=\fivesy
\textfont\itfam=\ninei \def\it{\fam\itfam\nineit}\def\sl{\fam\slfam\ninesl}%
\textfont\bffam=\ninebf \def\bf{\fam\bffam\ninebf}\rm}
%
%---------------------------------------------------------------------
%

\hyphenation{anom-aly anom-alies coun-ter-term coun-ter-terms}
\def\inv{^{\raise.15ex\hbox{${\scriptscriptstyle -}$}\kern-.05em 1}}

\def\Dsl{\,\raise.15ex\hbox{/}\mkern-13.5mu D} %this one can be subscripted
\def\dsl{\raise.15ex\hbox{/}\kern-.57em\partial}

\font\bigit=cmti10 scaled \magstep1
 %pound sterling
\def\lspace{\ifx\answ\bigans{}\else\qquad\fi}
\def\lbspace{\ifx\answ\bigans{}\else\hskip-.2in\fi} % $$\lbspace...$$
\def\boxeqn#1{\vcenter{\vbox{\hrule\hbox{\vrule\kern3pt\vbox{\kern3pt
    \hbox{${\displaystyle #1}$}\kern3pt}\kern3pt\vrule}\hrule}}}
\def\mbox#1#2{\vcenter{\hrule \hbox{\vrule height#2in
        \kern#1in \vrule} \hrule}}  %e.g. \mbox{.1}{.1}
%   matters of taste
%\def\tilde{\widetilde} \def\bar{\overline} \def\hat{\widehat}
%
% some sample definitions
  %     curly letters

\def\darr#1{\raise1.5ex\hbox{$\leftrightarrow$}\mkern-16.5mu #1}
 %pound sterling

 %puts a small half in a displayed eqn
\def\roughly#1{\raise.3ex\hbox{$#1$\kern-.75em\lower1ex\hbox{$\sim$}}}

\def\frac#1#2{{#1\over#2}}

\def\journal#1&#2(#3){\unskip, #1~\bf #2 \rm(19#3) }
\def\andjournal#1&#2(#3){\sl #1~\bf #2 \rm (19#3) }

\def\bra#1{\left\langle #1\right|}
\def\ket#1{\left| #1\right\rangle}

\catcode`\@=11\def\slash#1{\mathord{\mathpalette\c@ncel{#1}}}
\overfullrule=0pt
\def\steepslash{\c@ncel}
\def\frac#1#2{{#1\over #2}}

\def\:{\!:\!}
\def\inbar{\,\vrule height1.5ex width.4pt depth0pt}
\def\IQ{\relax\,\hbox{$\inbar\kern-.3em{\rm Q}$}}
\def\IB{\relax{\rm I\kern-.18em B}}
\def\IC{\relax\hbox{$\inbar\kern-.3em{\rm C}$}}
\def\IP{\relax{\rm I\kern-.18em P}}
\def\IR{\relax{\rm I\kern-.18em R}}
\def\ZZ{\relax\ifmmode\mathchoice
{\hbox{Z\kern-.4em Z}}{\hbox{Z\kern-.4em Z}}
{\lower.9pt\hbox{Z\kern-.4em Z}}
{\lower1.2pt\hbox{Z\kern-.4em Z}}\else{Z\kern-.4em Z}\fi}

\catcode`\@=12

%                      Zeitschriften:
\def\npb#1(#2)#3{{ Nucl. Phys. }{B#1} (#2) #3}
\def\plb#1(#2)#3{{ Phys. Lett. }{#1B} (#2) #3}
\def\pla#1(#2)#3{{ Phys. Lett. }{#1A} (#2) #3}
\def\prl#1(#2)#3{{ Phys. Rev. Lett. }{#1} (#2) #3}
\def\mpla#1(#2)#3{{ Mod. Phys. Lett. }{A#1} (#2) #3}
\def\ijmpa#1(#2)#3{{ Int. J. Mod. Phys. }{A#1} (#2) #3}
\def\cmp#1(#2)#3{{ Comm. Math. Phys. }{#1} (#2) #3}
\def\cqg#1(#2)#3{{ Class. Quantum Grav. }{#1} (#2) #3}
\def\jmp#1(#2)#3{{ J. Math. Phys. }{#1} (#2) #3}
\def\anp#1(#2)#3{{ Ann. Phys. }{#1} (#2) #3}
\def\prd#1(#2)#3{{ Phys. Rev. } {D{#1}} (#2) #3}
\def\ptp#1(#2)#3{{ Progr. Theor. Phys. }{#1} (#2) #3}
\def\aom#1(#2)#3{{ Ann. Math. }{#1} (#2) #3}

\def\br{\buildrel}
\def\bra{\langle}
\def\ket{\rangle}

\def\K{{\bf K}}

\def\P{{\bf P}}

\def\Z{{\bf Z}}

\def\cO{{\cal O}}

\def\cU{{\cal U}}

\def\cicy#1(#2|#3)#4{\left(\matrix{#2}\right|\!\!
                     \left|\matrix{#3}\right)^{{#4}}_{#1}}

\def\ra{\rightarrow}

\def\Box{{\,\lower0.9pt\vbox{\hrule
\hbox{\vrule height 0.2 cm \hskip 0.2 cm
\vrule height 0.2 cm}\hrule}\,}}

\global\newcount\thmno \global\thmno=0
\def\definition#1{\global\advance\thmno by1
\bigskip\noindent{\bf Definition \secsym\the\thmno. }{\it #1}
\par\nobreak\medskip\nobreak}
\def\question#1{\global\advance\thmno by1
\bigskip\noindent{\bf Question \secsym\the\thmno. }{\it #1}
\par\nobreak\medskip\nobreak}
\def\theorem#1{\global\advance\thmno by1
\bigskip\noindent{\bf Theorem \secsym\the\thmno. }{\it #1}
\par\nobreak\medskip\nobreak}
\def\proposition#1{\global\advance\thmno by1
\bigskip\noindent{\bf Proposition \secsym\the\thmno. }{\it #1}
\par\nobreak\medskip\nobreak}
\def\corollary#1{\global\advance\thmno by1
\bigskip\noindent{\bf Corollary \secsym\the\thmno. }{\it #1}
\par\nobreak\medskip\nobreak}
\def\lemma#1{\global\advance\thmno by1
\bigskip\noindent{\bf Lemma \secsym\the\thmno. }{\it #1}
\par\nobreak\medskip\nobreak}
\def\conjecture#1{\global\advance\thmno by1
\bigskip\noindent{\bf Conjecture \secsym\the\thmno. }{\it #1}
\par\nobreak\medskip\nobreak}
\def\exercise#1{\global\advance\thmno by1
\bigskip\noindent{\bf Exercise \secsym\the\thmno. }{\it #1}
\par\nobreak\medskip\nobreak}
\def\remark#1{\global\advance\thmno by1
\bigskip\noindent{\bf Remark \secsym\the\thmno. }{\it #1}
\par\nobreak\medskip\nobreak}
\def\problem#1{\global\advance\thmno by1
\bigskip\noindent{\bf Problem \secsym\the\thmno. }{\it #1}
\par\nobreak\medskip\nobreak}
\def\others#1#2{\global\advance\thmno by1
\bigskip\noindent{\bf #1 \secsym\the\thmno. }{\it #2}
\par\nobreak\medskip\nobreak}

\def\thmlab#1{\xdef #1{\secsym\the\thmno}\writedef{#1\leftbracket#1}\wrlabeL{#1=#1}}
%
% redefine \newsec so that all \thmno set to zero in a new section
%
\def\newsec#1{\global\advance\secno by1\message{(\the\secno. #1)}
%\ifx\answ\bigans \vfill\eject \else \bigbreak\bigskip \fi  %if desired
\global\subsecno=0\thmno=0\eqnres@t\noindent{\bf\the\secno. #1}
\writetoca{{\secsym} {#1}}\par\nobreak\medskip\nobreak}
\def\eqnres@t{\xdef\secsym{\the\secno.}\global\meqno=1\bigbreak\bigskip}
\def\sequentialequations{\def\eqnres@t{\bigbreak}}\xdef\secsym{}
%

%

%%%%%%%%%%%%%%%%%%%%%%
%  REFERENCES        %
%%%%%%%%%%%%%%%%%%%%%%

\ref\Atiyah{M. Atiyah, {\it Convexity and commuting Hamiltonians},
Bull. London Math. Soc. 14 (1982) 1-15.}
\ref\AB{M. Atiyah and R. Bott, {\it The moment map and
equivariant cohomology}, Topology 23 (1984) 1-28.}
\ref\Bertram{A. Bertram, {\it Another way to enumerate rational curves
with torus action}, math.AG/9905159.}
\ref\BehrendFentachi{K. Behrend and B. Fentachi,
{\it The intrinsic normal cone}, Invent. Math. 128 (1997) 45-88.}
\ref\BDPP{G. Bini, C. De Concini, M. Polito, and C. Procesi,
{\it Givental's work relative to mirror symmetry},
math.AG/9805097.}
\ref\Elezi{A. Elezi, {\it Mirror symmetry for concavex bundles
on projective spaces}, math.AG/0004157.}
\ref\FP{C. Faber, and R. Pandharipande, {\it Hodge Integrals and
Gromov-Witten Theory}, math.AG/9810173.}
\ref\Gathmann{A. Gathmann, {\it
Relative Gromov-Witten invariants and the mirror formula},
math.AG/0009190.}
%\ref\FriedmanMorgan{R. Friedman, and J.W. Morgan,
%{\it Smooth four-manifolds and complex surfaces,} in
%Ergebnisse der Mathematik und ihrer Grenzgebiete (3)
%[Results in Mathematics and Related Areas (3)], 27. Springer-Verlag, Berlin, 1994.}
\ref\Giv{A. Givental, {\it Equivariant Gromov-Witten
invariants}, alg-geom/9603021.}
\ref\GP{T. Graber and R. Pandharipande, {\it Localization of
virtual classes}, alg-geom/9708001.}
\ref\GS{V. Guillemin and S. Sternberg, {\it Convexity properties
of the moment mapping}, Invent. Math. 67 (1982) 491-513.}
\ref\Hirzebruch{F. Hirzebruch, {\it Topological methods in
algebraic geometry}, Springer-Verlag, Berlin 1995, 3rd Ed.}
\ref\GKM{M. Goresky, R. Kottwitz and R. MacPherson,
{\it Equivariant cohomology, Koszul duality and the
localization theorem}, Invent. Math. 131 (1998) 25-83.}
\ref\GP{T. Graber and R. Pandharipande, {\it Localization of
virtual classes}, alg-geom/9708001.}
\ref\KKV{S. Katz, A. Klemm, and C. Vafa, {\it Geometric engineering
of quantum field theories}, Nucl. Phys. B497 (1997) 173-195.}
%\ref\KM{F.F. Knudsen, and D. Mumford,
%{\it The projectivity of the moduli space of stable curves.
%I. Preliminaries on "det" and "Div". }
%Math. Scand. 39 (1976), no. 1, 19-55.}
\ref\K{M. Kontsevich,
{\it Enumeration of rational curves via torus actions.}
In: The Moduli Space of Curves, ed. by
R. Dijkgraaf, C. Faber, G. van der Geer, Progress in Math.
vol. 129, Birkh\"auser, 1995, 335--368.}
\ref\MP{D. Morrison and R. Plesser, {\it
Summing the instantons: quantum cohomology and mirror symmetry
in toric varieties}, alg-geom/9412236.}
\ref\LiTianII{ J. Li and G. Tian, {\it
Virtual moduli cycle and
Gromov-Witten invariants of algebraic varieties},
J. of Amer. math. Soc. 11, no. 1, (1998) 119-174.}
\ref\LLYI{B. Lian, K. Liu and S.T. Yau, {\it Mirror Principle I},
Asian J. Math. Vol. 1, No. 4 (1997) 729-763.}
\ref\LLYII{B. Lian, K. Liu and S.T. Yau, {\it Mirror Principle II},
Asian J. Math. Vol. 3, No. 1 (1999).}
\ref\LLYIII{B. Lian, K. Liu and S.T. Yau, {\it Mirror Principle III},
math.AG/9912038.}
\ref\LLYIV{B. Lian, K. Liu and S.T. Yau, {\it Mirror Principle IV},
math.AG/0007104.}
\ref\LCHY{B. Lian, C.H. Liu and S.T. Yau, {\it A Reconstruction of Euler Data},
math.AG/0003071.}
\ref\Pandha{R. Pandharipande, {\it Rational curves on
hypersurfaces (after givental)},
math.AG/9806133.}
\ref\Witten{E. Witten, {\it Phases of N=2 theories in two
dimension}, hep-th/9301042.}

\Title{}{A Survey of Mirror Principle}
 \centerline{\titlerms Bong H. Lian}
 \centerline{\it Department of Mathematics}                \vskip-1mm
 \centerline{\it Brandeis University, Waltham, MA 02154}   \vskip-1mm
 \centerline{ lian@brandeis.edu}         \vskip-1mm
\vskip .2in
\centerline{\titlerms Kefeng Liu}
 \centerline{\it Department of Mathematics}      \vskip-1mm
 \centerline{\it University of California, Los Angeles, CA 90024-6516}  \vskip-1mm
 \centerline{ kefeng@math.ucla.edu}                   \vskip-1mm
\vskip .2in
%\centerline{\titlerms and}
%\vskip .2in
 \centerline{\titlerms Shing-Tung Yau}
 \centerline{\it Department of Mathematics}                \vskip-1mm
 \centerline{\it Harvard University, Cambridge, MA 02138}  \vskip-1mm
 \centerline{ yau@math.harvard.edu}                  \vskip-1mm
\vfill
%\vskip .2in

Abstract.
This note briefly reviews the {\it Mirror Principle} as developed in the series of
papers \LLYI\LLYII\LLYIII\LLYIV\LCHY. We illustrate this theory with
a few new examples. One of them gives an intriguing connection to
a problem of counting holomorphic disks and annuli. This note has been submitted for
the proceedings of the Workshop on Strings, Duality and Geometry
the C.R.M. in Montreal of March 2000.

%\listtoc %input toc.tmp
%\writetoc %create toc.tmp

\Date{Sept 2000}

\newsec{Some Background}

In the aforementioned series of papers
we develop the {\it mirror principle} in increasing generality and breadth.
Given a projective manifold $X$, mirror principle is a theory
that yields relationships for and often computes the intersection
numbers of cohomology classes of the form
$b(V_D)$ on stable moduli spaces $\bar M_{g,k}(d,X)$.
Here $V_D$ is a certain induced vector bundles on $\bar M_{g,k}(d,X)$
and $b$ is any given multiplicative cohomology class.
In the first paper \LLYI, we consider this problem in the genus zero $g=0$ case
when $X=\P^n$ and $V_D$ is a bundle induced by any convex and/or concave bundle
$V$ on $\P^n$. As a consequence, we have proved a mirror formula which
computes the intersection numbers via a generating function.
When $X=\P^n$, $V$ is a direct sum of positive line bundles on $\P^n$,
and $b$ is the Euler class,
a second proof of this special case has been given
in \Pandha\BDPP~ following an approach proposed in \Giv.
Other proofs
in this case has also been given in \Bertram\Gathmann,
and when $V$ includes negative line bundles,
in \Elezi.
In \LLYII, we develop mirror principle
when $X$ is a projective manifold with $TX$ convex.
In \LLYIII, we consider the $g=0$ case when
$X$ is an arbitrary projective manifold.
Here  emphasis has been
put on a class of $T$-manifolds (which we call balloon manifolds)
because in this case mirror principle yields a (linear!)
reconstruction algorithm which
computes in principle all the intersection numbers above
for {\it any} convex/concave equivariant bundle $V$ on $X$ and {\it any}
equivariant multiplicative class $b$.
Moreover, specializing this theory to the case of
line bundles on toric manifolds and $b$ to Euler class,
we give a proof of the mirror formula for
toric manifolds.
In both \LLYIII~ and \LLYIV, we
develop mirror principle for higher genus.
We also extend the
theory to include the intersection numbers for
cohomology classes of the form $ev^*(\phi)b(V_D)$.
Here $ev:\bar M_{g,k}(d,X)\ra X^k$ is the usual evaluation map
into the product $X^k$ of $k$ copies of $X$, and $\phi$ is any cohomology
class on $X^k$.

For motivations and some historical background of the mirror principle,
we refer the reader to the introduction of \LLYI\LLYII.

In section 2, we outline the main ideas of the mirror principle,
and explain one of our main theorems.
In section 3, we discuss a few examples.

{\it Acknowledgment:} We thank C. Vafa for
informing us of his result on the disk counting problem. We also thank
the organizers for inviting us to
the conference on Geometry and String Theory
at the C.R.M. in Montreal in 2000.
B.H.L.'s research
is supported by NSF grant DMS-0072158.
K.L.'s research is supported by NSF grant
DMS-9803234 and the Terman fellowship and the Sloan fellowship.
S.T.Y.'s research is supported by DOE grant
DE-FG02-88ER25065 and NSF grant DMS-9803347.

\newsec{Mirror Principle}

For simplicity, we restrict our discussions to the genus zero theory,
and refer the interested reader to \LLYIV~ for a theory of higher genus.
Let $X$ be a projective $n$-fold, and $d\in H^2(X,\Z)$.
Let $M_{0,k}(d,X)$ denote the moduli stack of
$k$-pointed, genus 0, degree $d$, stable maps
$(C,f,x_1,..,x_k)$ on $X$ \K.
(Note that our notation is without the bar.)
By \LiTianII~(cf. \BehrendFentachi),
each nonempty
$M_{0,k}(d,X)$ admits a homology cycle $LT_{0,k}(d,X)$
of degree $dim~X+\bra c_1(X),d\ket +k-3$. This
cycle plays the role of the fundamental class in
topology, hence
$LT_{0,k}(d,X)$ is
called the virtual fundamental class.

Let $V$ be a convex vector bundle on $X$.
(ie. $H^1(\P^1,f^*V)=0$ for every holomorphic
map $f:\P^1\ra X$.) Then $V$ induces on each
$M_{0,k}(d,X)$ a vector bundle $V_d$,
with fiber at
$(C,f,x_1,..,x_k)$ given by the section space $H^0(C,f^*V)$.
Let $b$ be any multiplicative characteristic class
\Hirzebruch.
(ie. if $0\ra E'\ra E\ra E''\ra 0$ is an exact sequence
of vector bundles,
then $b(E)=b(E')b(E'')$.)
The problem we study here is to compute
the characteristic numbers
$$K_d:=\int_{LT_{0,0}(d,X)} b(V_d)$$
and their generating function:
$$\Phi(t):=\sum K_d~ e^{d\cdot t}.$$

There is a similar and equally important problem if one
starts from a concave vector bundle $V$ \LLYI.
(ie. $H^0(\P^1,f^*V)=0$ for every holomorphic
map $f:\P^1\ra X$.) More generally, $V$ can
be a direct sum of a convex and a concave bundle.

The rough idea of the Mirror Principle is that
the classes the induced bundles $V_d$ on the stable
moduli inherit a number of universal structures (ie.
exist in all stable map moduli of any projective manifold).
These structures combined with the multiplicative
properties of the classes $b(V_d)$ give rise to
some remarkable quadratic identities. It is often the case (when sufficient
symmetry is present on $X$) that
these identities are strong enough for a complete reconstruction of
the intersection numbers $K_d$.
We explain this idea further below without proofs.
For details, see \LLYIII.

{\it Step 1. Localization on the linear sigma model.}
Consider the moduli spaces
$M_d(X):=M_{0,0}((1,d),\P^1\times X)$. The projection
$\P^1\times X\ra X$ induces a map $\pi:M_d(X)\ra M_{0,0}(d,X)$.
Moreover, the standard action of $S^1$ on $\P^1$ induces
an $S^1$ action on $M_d(X)$. We first study a
slightly different problem. Namely
consider the classes $\pi^*b(V_d)$ on $M_d(X)$,
instead of $b(V_d)$ on $M_{0,0}(d,X)$.
First, there
is a canonical way to
embed fiber products
$$F_r=M_{0,1}(r,X)\times_X M_{0,1}(d-r,X)$$
each as an $S^1$ fixed point component into $M_d(X)$.
Let $i_r:F_r\ra M_d(X)$ be the inclusion map.
Second, there is an evaluation map
$e:F_r\ra X$ for each $r$.
Third, there is a (product of) projective space $W_d$ equipped with an $S^1$ action,
and there is an equivariant map $\varphi:M_d(X)\ra W_d$, and
embeddings $j_r:X\ra W_d$, such that the diagram
$$\matrix{
F_r & {\br i_r\over \longrightarrow} & M_d(X)\cr
e\downarrow &  & \downarrow \varphi\cr
X & {\br j_r\over \longrightarrow} & W_d}
$$
commutes.
Let $\alpha$ denotes the weight of
the standard $S^1$ action on $\P^1$.
Applying the localization formula \AB\GP, this diagram allows us
to recast our problem to one of studying
the $S^1$-equivariant classes
$$Q_d:=\varphi_*\pi^*b(V_d)$$
defined on $W_d$.
%Moreover we can expand the class
%$$A_d:={j^*_0Q_d\over e_{S^1}(X_0/W_d)}$$
% on $X$
%in powers of $\alpha^{-1}$,
%and find that it is of order $\alpha^{-2}$.

The projective space $W_d$ in the commutative diagram above
is called a linear sigma model of $X$.
They have been introduced in \MP~following \Witten.

{\it Step 2. Gluing identity.}
Consider the vector bundle $\cU_d:=\pi^* V_d\ra M_d(X)$,
restricted to the fixed point components $F_r$.
A point in $(C,f)$ in $F_r$ is
a pair $(C_1,f_1,x_1)\times(C_2,f_2,x_2)$ of
1-pointed stable maps glued together at
the marked points, ie. $f_1(x_1)=f_2(x_2)$.
From this, we get an exact sequence of bundles
on $F_r$:
$$0\rightarrow i_r^*\cU_d\rightarrow U_r'\oplus U_{d-r}'\rightarrow
e^*V\rightarrow 0.$$
Here $i_r^*\cU_d$ is the restriction to $F_r$ of the bundle
 $\cU_d\ra M_d(X)$. And $U_r'$ is the pullback of the bundle
$U_r\ra M_{0,1}(d,X)$ induced by $V$, and similarly for $U_{d-r}'$.
Taking the multiplicative characteristic class $b$, we get
the identity on $F_r$:
$$e^*b(V)b(i^*_r\cU_d)=b(U_r')b(U_{d-r}').$$
This is what we call
the {\it gluing identity}.
This may be translated to a similar quadratic identity, via
Step 1, for $Q_d$ in the equivariant Chow
groups of $W_d$.
The new identity is called the Euler data identity.

{\it Step 3. Linking theorem.}
The construction above is functorial, so that if $X$
comes equipped with a torus $T$ action, then the entire
construction becomes $G=S^1\times T$ equivariant and
not just $S^1$ equivariant. In particular,
the Euler data identity is an identity of
$G$-equivariant classes on $W_d$.
Our problem is  to first
compute the $G$-equivariant classes $Q_d$ on $W_d$
satisfying the Euler data identity.
Note that  the restrictions $Q_d|_p$ to the $T$ fixed points
$p$ in $X_0\subset W_d$ are polynomials functions on
the Lie algebra of $G$.
Suppose that $X$ is a balloon manifold.
This is a complex projective $T$-manifold
satisfying the following conditions \GKM:
\item{(i)} The $T$ fixed points are isolated.
\item{(ii)} Let $p$ be a $T$ fixed point. Then
the $T$ weights
$\lambda_1,..,\lambda_n$ of
the isotropic representation on the tangent space $T_pX$
are pairwise linearly independent.
We further assume that the moment map is 1-1 on the fixed point set.

In this case, the classes $Q_d$ are uniquely determined by the
values of the $Q_d|_p$, when $\alpha$ is some scalar multiple of a
weight $\lambda_i$. These values of $Q_d|_p$, which we call the
linking values (see \LLYIII~ for precise definition), can be
computed explicitly by exploiting the moment map \Atiyah\GS~ as
well as certain structure of a balloon manifold.

\theorem{\LLYIII The equivariant classes $Q_d=\varphi_*\pi^*b(V_d)$,
as a solution to the Euler data identity,
can be completely recovered from the linking values.}

Once the linking values are known, it is often easy to manufacture
explicitly the $G$-equivariant classes $Q_d$ using the linking
values as a guide. Many explicit examples are discussed in \LLYIII.

{\it Step 4. Computing the $K_d$.}
Once the classes $Q_d=\varphi_*\pi^* b(V_d)$ are determined,
one has to carefully unwind the the commutative diagram and maps in Step 1.
This allows us to establish a crucial integral identity
between $K_d$ and the classes $Q_d$, which in turn allows us
to compute the $K_d$.

\newsec{Some Examples}

We now give some examples of our main theorem.

{\it Line bundles on a balloon manifold.}
Let $X$ be balloon manifold, and set $b$ to be the Chern polynomial.
We fix a base $H_1,..,H_m$ of $H^2(X)$, and let $t_1,..,t_m$ be formal variables.
We denote by $e_{S^1}(X_0/W_d)$ the equivariant Euler class of the normal
bundle of $X_0:=j_0(X)\subset W_d$.
Let
$$V=V^+\oplus V^-,~~~V^+:=\oplus L_i^+,~~~V^-:= \oplus L_j^-$$
satisfying
$c_1(V^+)-c_1(V^-)=c_1(X)$, and $rk~V^+-rk~V^- -(n-3)\geq0$.
where the $L_i^\pm$ are respectively convex/concave line
bundles on $X$. Let
$$\eqalign{
\Omega&=B_0:=c(V^+)/c(V^-)=\prod_i(x+c_1(L^+_i))
/\prod_j(x+c_1(L^-_j))\cr
B_d&:={1\over e_{S^1}(X_0/W_d)}\times
\prod_i\prod_{k=0}^{\bra c_1(L^+_i),d\ket}
(x+c_1( L^+_i)-k\alpha)\times
\prod_j\prod_{k=1}^{-\bra c_1(L^-_j),d\ket-1}
(x+c_1( L^-_j)+k\alpha).\cr
B(t)&:=e^{-H\cdot t}\sum B_d e^{d\cdot t}\cr
\Phi(t)&:=\sum K_d e^{d\cdot t}.
}
$$

\theorem{There exist unique power series $f(t),g(t)$ such that
the following formula holds:
$$
{1\over s!}\left({d\over dx}\right)^s|_{x=0}
\int_X\left( e^{f/\alpha} B(t)-
e^{-H\cdot\tilde t /\alpha}\Omega\right)
=\alpha^{-3}x^{-s}(2\Phi(\tilde t)
-\sum_i\tilde t_i{\partial\Phi(\tilde t)\over\partial\tilde t_i}).
$$
where $s:=rk~V^+-rk~V^--(n-3)$, $\tilde t:=t+g$.
Moreover, $f,g$ are determined by the condition that the
integrand on the left hand side is of order $O(\alpha^{-2})$.
}

Note that when $x\ra 0$, the formula above reduces to the
case when $b$ is the Euler class.

{\it The tangent bundle on $\P^n$.}
The example above deals, of course, with direct sum of line bundles
only. We now give an example starting from the tangent bundle $V=TX$
on $X=\P^n$, which is nonsplit. Consider the case
where $b_T$ the $T$-equivariant Chern polynomial.
Let $\lambda_i$ be the weights of the standard $T$ action on $\P^n$,
and $p_i$ be the $i$th fixed point.
Recall that $\Omega:=b_T(V)={1\over x}\prod_i(x+H-\lambda_i)$, where
$H$ is the equivariant class on $\P^n$ with $H|_{p_i}=\lambda_i$.

Using the $T$ equivariant Euler sequence
$$
0\ra\cO\ra\oplus_{i=0}^n\cO(H-\lambda_i)\ra TX\ra0
$$
One can compute the linking values in this case. There are given by
$$
\prod_i\prod_{k=0}^d(x+\lambda_j-\lambda_i-k\lambda/d).
$$
Here $p,q$ are the $j$th and the $l$th fixed points in
$\P^n$, and $\lambda=\lambda_j-\lambda_l$.
We can use this to set up a system of linear equations to solve for $A(t)$
inductively. However, there is an easier way to compute $A(t)$
in this case. Using the linking values above as a guide, we set
$$B_d:={1\over x}\prod_i\prod_{k=0}^d(x+H-\lambda_i-k\alpha)$$
and let
$$B(t):=e^{-H\cdot t/\alpha}\sum {B_d\over \prod_i\prod_{k=1}^d(H-\lambda_i-k\alpha)}.$$
Then it can be shown that the series
$$
A(t):=e^{-H\cdot t/\alpha}\sum {j_0^*Q_d\over \prod_i\prod_{k=1}^d(H-\lambda_i-k\alpha)}
 e^{d\cdot t}
$$
is related to $B(t)$ by
$$A(t+g)=e^{f/\alpha} B(t)$$
where $f,g$ are explicitly computable functions,
similar to those in the previous example.
This relation, once again, allows us to compute all the $K_d$
simultaneously.

{\it $V:=\cO(-1)\oplus\cO(-1)$ on $\P^1$.} In this case, we let $b$ be
the Euler class $c_{top}$, and we would like to compute
the one-pointed intersection numbers
$$\int_{M_{0,1}(d,\P^1)} e^*(H)b(V_d').$$
Here $V_d'$ is the bundle induced on $M_{0,1}(d,\P^1)$ by $V$,
$H$ is the hyperplane class on $\P^1$, and $e:M_{0,1}(d,\P^1)\ra\P^1$
is the evaluation map. We can easily specialize the first example above to the
case of $X=\P^1$ and $V=\cO(-1)\oplus\cO(-1)$. In this case, $f=g=0$,
and we get the formula
$$\int_{\P^1} e^{-Ht/\alpha} {j^*_0 Q_d\over \prod_{m=1}^d(H-m\alpha)^2}
=\alpha^{-3}(2-dt)K_d.$$
Apply ${d\over dt}$ to both sides, and combine the result with
Theorem 3.2 in \LLYI ~(see the first eqn. on p36 there). What we get is
$$\int_{M_{0,1}(d,\P^1)} e^*(H)b(V_d')=d~K_d=d^{-2}.$$

The values $d^{-2}$ remind us of a result that Vafa obtains via
the physics of local mirror symmetry. He considers the problem of
counting holomorphic disks in a Calabi-Yau 3-fold equipped with a
choice of Lagrangian submanifold. The boundary of the disks are
required to lie in the Lagrangian submanifold. A counting problem
is heuristically formulated into a problem of determining the
Euler classes of certain yet-to-be-defined moduli spaces. In this
case the physics of local mirror symmetry indicates that the Euler
classes should be a given by
a ``multiple-cover'' contribution $d^{-2}$, where $d$ is the winding number of the
disk's boundary circle along the Lagrangian submanifold.

Vafa's result suggests the following interpretation. 
The Lagrangian submanifold plays the role of a vanishing
cycle in a certain limit. A holomorphic disk with boundary
landing on the Lagrangian submanifold would look like
a $\P^1$ with one marked point in this limit.
 Local mirror symmetry
suggests that we should use the stable map moduli spaces of $\P^1$
as a model for this problem. The requirement that the marked point
lands on the vanishing cycle may be thought
of as the incidence condition on the map $\P^1\ra\P^1$ with one
point mapped to the cycle $H$. The appropriate
moduli spaces in this model should then be $M_{0,1}(d,\P^1)$, and
the Euler classes should correspond to the ``multiple-cover'' formula for the
bundle $V_d'$ induced by $V=\cO(-1)\oplus\cO(-1)$. So a good
candidate for the intersection numbers are
 $\int_{M_{0,1}(d,\P^1)}
e^*(H)b(V_d')=d^{-2}$.

{\it Two Lagrangian $S^3$ in a CY 3-fold?.} Another interesting
situation considered by physicists is the problem of counting
annuli in a CY 3-fold equipped with two Lagrangian 3-spheres $S^3$,
subject to the incidence condition that each of the
boundary circles of the annulus lands inside one of the 3-spheres.
Again the 3-spheres plays the role of two vanishing cycle,
and are allowed to contract
to points $x,y$. The annulus looks like a $\P^1$ with two marked points
anchored to $x,y$. By analogy with the previous example as
in local mirror symmetry,
the corresponding stable map moduli in this case should be
the two-pointed moduli $M_{0,2}(d,\P^1)$, and the corresponding
intersection numbers should be
$$\int_{M_{0,2}(d,\P^1)} e_1^*(H)e_2^*(H)b(V_d'').$$
Here $V_d''$ is the induced bundle $\rho_2^*V_d'$ on $M_{0,2}(d,\P^1)$,
where $\rho_2: M_{0,2}(d,\P^1)\ra M_{0,1}(d,\P^1)$ is the
map that forgets the second marked point, and the $e_i$
are the usual evaluation maps on $M_{0,2}(d,\P^1)$.

The intersection number can be easily computed in a way
analogous to the previous example. By writing $b(V_d'')=\rho_2^* b(V_d')$,
we get
$$\int_{M_{0,2}(d,\P^1)} e_1^*(H)e_2^*(H)b(V_d'')=
\int_{M_{0,1}(d,\P^1)} e^*(H)b(V_d'){\rho_2}_* e_2^*(H).$$
By integrating along a fiber of the map
$\rho_2$, we see that the last factor in the integrand
contributes an overall factor $d$.
Thus we get the answer
$$\int_{M_{0,2}(d,\P^1)} e_1^*(H)e_2^*(H)b(V_d'')=d^{-1}.$$

The last two example suggests the very interesting possibility
that one may be able to use stable map moduli spaces as models for
some of the moduli spaces in the problem of counting holomorphic
disks and annuli with suitable incidence conditions. Moreover, the
appropriate intersection numbers should come from Euler classes of
induced bundles, which is exactly what the mirror principle is designed
to study. This possibility deserves
further investigations.

\footatend\vfill\supereject\immediate\closeout\rfile\writestoppt
\baselineskip=14pt\centerline{{\bf References}}\bigskip{\frenchspacing%
\parindent=20pt\escapechar=` \input refs.tmp\vfill\eject}\nonfrenchspacing

\end